\documentclass[11pt]{amsart}
\usepackage{verbatim, graphicx, rotating}

\usepackage {amssymb}

\input xy
\xyoption{all}
\input epsf

\newlength{\tabwidth}
\newlength{\tabheight}
\newlength{\tabrule}
\newlength{\tabwidthx}
\newlength{\tabheightx}

\def\gentabbox#1#2#3#4{\vbox to \tabheight{\setlength{\tabrule}{#3}%
  \setlength{\tabwidthx}{#1\tabwidth}\addtolength{\tabwidthx}{\tabrule}%

\setlength{\tabheightx}{#2\tabheight}\addtolength{\tabheightx}{-\tabheight}%
  \hbox to #1\tabwidth{%
    \hspace{-0.5\tabrule}\rule{\tabrule}{#2\tabheight}\hspace{-\tabrule}%
    \vbox to #2\tabheight{\hsize=\tabwidthx%
      \vspace{-0.5\tabrule}\hrule width\tabwidthx height\tabrule%
      \vspace{-0.5\tabrule}\vfil%
      \hbox to \tabwidthx{\hss#4\hss}%
        \vfil\vspace{-0.5\tabrule}%
      \hrule width\tabwidthx height\tabrule\vspace{-0.5\tabrule}}%
    \hspace{-\tabrule}\rule{\tabrule}{#2\tabheight}\hspace{-0.5\tabrule}}%
  \vspace{-\tabheightx}}}
\def\genblankbox#1#2{\vbox to \tabheight{\vfil\hbox to
#1\tabwidth{\hfil}}}
\def\tabbox#1#2#3{\gentabbox{#1}{#2}{0.4pt}{\strut #3}}

\catcode`\:=13 \catcode`\.=13 \catcode`\;=13 
\catcode`\>=13 \catcode`\^=13
\def:#1\\{\hbox{$#1$}}
\def.#1{\tabbox{1}{1}{$#1$}}
\def>#1{\tabbox{2}{1}{$#1$}}
\def^#1{\tabbox{1}{2}{$#1$}}
\def;{\genblankbox{1}{1}\relax}
\catcode`\:=12 \catcode`\.=12 \catcode`\;=12 
\catcode`\>=12 \catcode`\^=7

\newcommand{\iIm}{\text{A}}
\newcommand{\jIm}{\text{B}}

\newcommand{\cl}{{cl}}

\newcommand{\ol}{\overline}

\newcommand{\field}{\mathbb}
\newcommand{\liealgebra}{\mathfrak}
\newcommand{\la}{\liealgebra}

\newcommand{\C}{{\field C}}
\newcommand{\R}{{\field R}}

\renewcommand{\b}{\liealgebra b}

\newcommand{\n}{{\la n}}

\newcommand{\lra}{\longrightarrow}

\newcommand{\bs}{\backslash}

\newcommand{\Sp}{\mathrm{Sp}}
\newcommand{\PSp}{\mathrm{PSp}}
\newcommand{\Spin}{\mathrm{Spin}}
\newcommand{\SO}{\mathrm{SO}}
\newcommand{\SL}{\mathrm{SL}}
\renewcommand{\O}{\mathrm{O}}
\newcommand{\U}{\mathrm{U}}
\newcommand{\Spq}{\Sigma_\pm(p,q)}
\newcommand{\SPpq}{\Sigma^\sym_\pm(2p,2q)}
\newcommand{\SOpq}{\Sigma^\asym_\pm(n,n)}
\newcommand{\SOnn}{\Sigma^\asym_\pm(n,n)}
\newcommand{\sym}{\mathrm{sym}}
\newcommand{\ad}{\mathrm{ad}}
\renewcommand{\sc}{\mathrm{sc}}
\newcommand{\asym}{\mathrm{asym}}
\newcommand{\SU}{\mathrm{SU}}

\newcommand{\GL}{\mathrm{GL}}

\newcommand{\+}{\! + \!}

\newtheorem{prop}{Proposition}[section]

\newtheorem{lemma}[prop]{Lemma}
\newtheorem{theorem}[prop]{Theorem}

\theoremstyle{definition}
\newtheorem{defn}[prop]{Definition}
\newtheorem{remark}[prop]{Remark}
\newtheorem{example}[prop]{Example}

\newcommand{\frb}{\mathfrak{b}}

\newcommand{\frg}{\mathfrak{g}}
\newcommand{\frh}{\mathfrak{h}}

\newcommand{\frl}{\mathfrak{l}}
\newcommand{\frn}{\mathfrak{n}}
\newcommand{\fro}{\mathfrak{o}}
\newcommand{\frp}{\mathfrak{p}}

\newcommand{\frs}{\mathfrak{s}}

\newcommand{\bbC}{\mathbb{C}}

\newcommand{\bbH}{\mathbb{H}}

\newcommand{\bbN}{\mathbb{N}}

\newcommand{\bbR}{\mathbb{R}}

\newcommand{\bbZ}{\mathbb{Z}}

\newcommand{\caB}{\mathcal{B}}

\newcommand{\caO}{\mathcal{O}}
\newcommand{\caP}{\mathcal{P}}

\begin{document}
\title[pattern avoidance and smoothness]
{Pattern avoidance and smoothness of
closures for orbits of a symmetric subgroup 
in the flag variety}

\author{William M.~McGovern and Peter E.~Trapa}
\date{\today}
\thanks{The second author was partially supported by
NSA grant MSPF-06Y-096 and NSF grant DMS-0532393.}

\address{Department of Mathematics, University of Washington, Seattle, WA
98195}

\email{mcgovern@math.washington.edu}

\address{Department of Mathematics, University of Utah, Salt Lake City, UT
84112}

\email{ptrapa@math.utah.edu}

\begin{abstract}
We give a pattern avoidance criterion to 
classify the orbits of $\Sp(p,\C) \times \Sp(q,\C)$ 
(resp.~$\GL(n,\C)$) on the flag variety of type $\mathrm{C}_{p\+q}$
(resp.~$\mathrm{D}_n$) with rationally smooth closure.
We show that all such orbit closures fiber (with smooth fiber) over a 
smaller flag variety, and hence are in fact smooth.  In addition we
prove that the classification is insensitive to isogeny.
\end{abstract}

\maketitle

Suppose $G$ is a complex connected reductive algebraic group
and let $\theta$ denote an involutive automorphism of $G$.
Write $K$ for the fixed points of $\theta$, and $\caB$ for variety of
maximal solvable subalgebras of the Lie algebra $\frg$ of $G$.  (Henceforth
we call this variety simply the flag variety of $G$.)  Then $K$ acts
with finitely many orbits on $\caB$ via the restriction of the adjoint
action (e.g.~\cite{matsuki}).

Since we have assumed the ground field is $\bbC$,
$\theta$ arises as the
complexification of a Cartan involution for a real form $G_\R$ of
$G$.  The localization theory of Beilinson-Bernstein relates the
geometry of $K$ orbits on $\caB$ with the category of Harish-Chandra
modules for $G_\R$.  Meanwhile, as a  special case, one can
consider the setting where $G_\R$ is itself a complex Lie group.  In
this case $G = G_\R \times G_\R$, $\theta$ is the involution that
interchanges the two factors, $K$ is the diagonal copy of $G_\R$, and
the Weyl group of $G_\R$ parametrizes the orbits of $K$ on $\caB$
(which now is two copies of the flag variety $\caB_\circ$ for $G_\R$).
Intersecting such an orbit with (say) the left copy
$\caB_\circ$ gives an orbit of a Borel subgroup on $\caB_\circ$, that
is a Schubert cell.  This process preserves the fine structure of the
singularities of the closures of each kind of orbit, and is the
geometric underpinning of the 
equivalence of
categories (essentially) between category $\caO$ and a
suitable category of Harish-Chandra modules for $G_\R$ (e.g.~\cite{BB}).

Thus, roughly speaking, any question which one can ask about Schubert
varieties can also be posed in the greater generality of $K$ orbits on
$\caB$, and 
any relation of the geometry of the former with category
$\caO$ can potentially be translated into a relation of the latter
with the category of Harish-Chandra modules.  
Some of the deepest results in this direction are due to Lusztig-Vogan
\cite{LV} and Vogan \cite{vogan:ic2},\cite{vogan:ic3}, which when
taken together give an algorithm to compute the local intersection
homology (with coefficients in any irreducible local system) of any
orbit closure.  When $G_\R$ is a complex group, the algorithm is
equivalent to that of \cite{KL} for Schubert varieties.

Our interest here is determining when the closures of $K$ orbits on
$\caB$ are smooth, or more generally rationally smooth.  Since the latter
condition is equivalent to a condition on local intersection homology
(see the discussion of \cite[Appendix]{KL}), the question of whether a particular orbit has rationally
smooth closure can be answered using the algorithm of \cite{vogan:ic2}.
But it is desirable to have a ``closed form'' of the answer.  For
instance, in the case that $G_\R$ is complex (or, equivalently, the
case of Schubert varieties), a closed form answer for smoothness
or rational smoothness has been obtained in terms of a kind of 
pattern avoidance
for Weyl group elements\footnote{There is an extensive literature of geometric and combinatorial
results relating
pattern avoidance in Weyl groups to singularities of Schubert varieties.  We
do not attempt to recount the history of these results in detail here.  See  \cite[Chapter 8]{BL} 
or \cite{BP} and the extensive references therein.}.

In the setting of $K$ orbits on $\caB$, new phenomena appear that are not present
in the case of Schubert varieties.  The most obvious difference is that while the definition of Schubert varieties
is independent of the isogeny class of $G$, symmetric subgroups (and their orbits on $\caB$)
do indeed depend on isogeny.   A typical complication in the latter case may be
visualized as follows.  Suppose $G$ is simply connected but not adjoint, $\ol G\neq G$ is isogenous to $G$,
and $\ol K$ and $K$ are respective symmetric subgroups with the same Lie algebra.  Then it frequently
happens that $\ol K$ is disconnected (while $K$ is connected)\footnote{We remark that the notion
of isogeny we are considering here differs from isogeny of the corresponding real forms.  More precisely,
if $G_\R$ and $\ol G_\R$ are the real forms corresponding to symmetric subgroups $K$ and $\ol K$ of groups $G$
and $\ol G$ with $G$ simply connected and $\ol G$ a quotient of $G$ by a central subgroup, then of course it does not follow that
$\ol G_\R$ is a quotient of $G_\R$ by a central subgroup.  (In other words, the relevant central subgroup of $G$ need
not be defined over $\bbR$.) 
For the results of this paper, the latter notion of isogeny of real forms is not interesting and will not be considered.}.  
Thus there may be two distinct orbits 
$Q_1$ and $Q_2$ for $K$ on $\caB$ whose union forms a single
orbit $Q$ for $\ol K$.   Schematically one encounters pictures as follows.
 
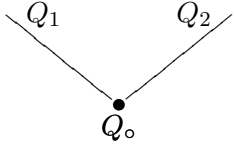
\[
\xy
(0,10)*{Q_1}+(-5,0); (10,-5)*{Q_\circ}+(0,3)*{\bullet} **\dir{-};
(10,-5)*{Q_\circ}+(0,3)*{\bullet}; (20,10)*{Q_2}+(5,0)**\dir{-};
\endxy
\]
Here $Q_\circ$ is a closed $K$ orbit (which is also a $\ol K$ orbit) which appears in the closure of both $K$
orbits $Q_1$ and $Q_2$.  The picture indicates that the closure of the $K$ orbit $Q_1$ (or $Q_2$) is smooth at $Q_\circ$.  But the 
$\ol K$ orbit $Q$, the union of $Q_1$ and $Q_2$, has closure which isn't smooth (or even rationally smooth) at $Q_\circ$.
Moreover, this particular example can ``propagate'' in higher rank schematically as follows.
 
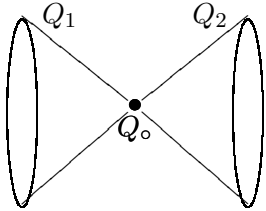
\[
 \xy
(-5,-3)*\xycircle(2,12.5){-};
(25,-3)*\xycircle(2,12.5){-};
(0,10)*{Q_1}+(-5,0); (10,-5)*{Q_\circ}+(0,3)*{\bullet} **\dir{-};
(10,-5)*{Q_\circ}+(0,3)*{\bullet}; (20,10)*{Q_2}+(5,0)**\dir{-};
(10,-5)+(0,3)*{\bullet}; (-5,-15)**\dir{-};
(10,-5)+(0,3)*{\bullet}; (25,-15)**\dir{-};
\endxy
\]
This time the closure of the $K$ orbit $Q_1$ (or $Q_2$) is no longer smooth at $Q_\circ$ (but
is rationally smooth), but once again the 
$\ol K$ orbit $Q$ has closure which is not smooth (nor rationally smooth) at $Q_\circ$. 
Of course none of this is conceptually very complicated, but
since, roughly speaking, pattern avoidance
results are predicated on failure of smoothness propagating uniformly to higher rank,
these examples suggest that such results are potentially more 
subtle in the case of $K$ orbits
in ways that are not seen for Schubert varieties.   
This is made clear by Example \ref{e:bad}.

In this paper, we are interesting in understanding nice isogeny-independent pattern 
avoidance results.
For instance, \cite{mcg:upq} answers
the questions of smoothness and rational smoothness for the closures
of orbits of $K=\GL(p,\C) \times \GL(q,\C)$ on the flag variety for
$G= \GL(p\+q,\C)$.  This is the setting arising from the real form
$G_\R = \U(p,q)$ of $G$.  
After reviewing some preliminaries in Section \ref{s:prelim},
we recall the results of \cite{mcg:upq} in Section \ref{s:upq}
and then prove that they are,
in a suitable sense, insensitive to isogeny.

In the remainder of the paper, we go on to study other classical groups outside of type A.
In general, as remarked above, the situation potentially depends crucially on issues
related to isogeny.
But, perhaps surprisingly, we find that isogeny essentially plays no role
if $G$ is of type C (resp.~type D) and the Lie algebra of $K$ is
$\frs\frp(p,\C) \oplus \frs\frp(q,\C)$ (resp.~$\frg\frl(n,\C)$).  This setting includes the
cases arising from the real groups $\Sp(p,q)$ and $\SO^*(2n)$.
The main results are Theorems \ref{t:sppq} and \ref{t:sostar}.  They
completely settle the question of rational smoothness (which, in the end, turns out
to be equivalent to smoothness) for $K$ orbits closures in these cases.
The statements are formulated in terms of a remarkably 
simple pattern avoidance criterion
resembling the one discovered in \cite{mcg:upq}.  (Roughly speaking there are
always only eight ``bad'' patterns to avoid.)

\section{preliminaries}
\label{s:prelim}

In particular examples below, we will need a detailed description of
the closure order of $K$ orbits on $\caB$.  We begin by
recalling a few features of the general case.  These are due
to Matsuki \cite{matsuki}, Matsuki-Oshima \cite{MO}, Lusztig-Vogan
\cite{LV}, and
are given a full exposition in Richardson-Springer \cite{RS}.

Following the terminology of \cite[Section 5.1]{RS}, we first recall the
weak closure order on $K \bs \caB$.  Fix, once and for all, a choice
of $\theta$-stable Cartan subalgebra $\frh$ in $\frg$,
and a choice of positive roots $\Delta^+$ in the full root system
$\Delta = \Delta(\frg,\frh)$.  
Write $B$ for the corresponding Borel subgroup of $G$.
For a simple root $\alpha \in \Delta^+$, let
$\caP_\alpha$ denote the variety of parabolic subalgebras of type $\alpha$
 (i.e.~those conjugate to the one with roots
$\Delta^+\cup\{-\alpha\}$).  Write $\pi_\alpha$ for the projection of
$\caB$ to $\caP_\alpha$, and define
\begin{equation}
\label{e:pnot}
s_\alpha \cdot \caO := \text{ the (unique) dense $K$ orbit in } 
\pi_\alpha^{-1}(\pi_\alpha(\caO)).
\end{equation}
The weak closure order is generated by relations
$\caO < s_\alpha\cdot \caO$ whenever $\dim(s_\alpha\cdot \caO) = 
\dim(\caO) + 1$.  In this case, for $\caO' = s_\alpha \cdot \caO$, we write
\[
\caO \stackrel{\alpha}{\lra} \caO'.
\]
Then $\caO \subset {\ol \caO'}$, but these relations do not generate
the full closure order.  To obtain all closure relations, we must
recursively apply the following procedure (implicit in
\cite[Theorem 7.11(vii)]{RS}, for instance).
Whenever a codimension one subdiagram
of the form
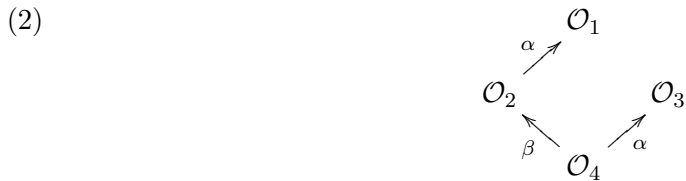
\begin{equation}
\label{e:exchange}
{
\xymatrixcolsep{1pc} 
\xymatrixrowsep{1pc}
\xymatrix
{& \caO_1 \\
\caO_2 \ar@{>}[ur]^\alpha 
& & \caO_3\\
& \caO_4 \ar@{>}[ul]^\beta \ar@{>}[ur]_\alpha
}
}
\end{equation}
is encountered, it must be completed to
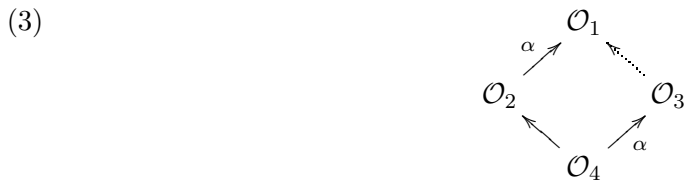
\begin{equation}
\label{e:exchange2}
{
\xymatrixcolsep{1pc} 
\xymatrixrowsep{1pc}
\xymatrix
{& \caO_1 \\
\caO_2 \ar@{>}[ur]^\alpha 
& & \caO_3\ar@{.>}[ul]\\
& \caO_4 \ar@{>}[ul] \ar@{>}[ur]_\alpha
}
}
\end{equation}
New edges added in this way are represented by dashed lines
in the figures appearing in Section \ref{s:examples}.
Note that this operation must be applied recursively,
and thus the solid unlabeled edge in the original diagram \eqref{e:exchange}
may be dashed as the recursion unfolds.

We next recall a definition from \cite[4.7]{RS}.  Given
$w \in W = W(\frh,\frg)$, fix a reduced expression $w = s_{\alpha_k}\cdots s_{\alpha_2}s_{\alpha_1}$.
For $\caO \in K\bs \caB$, define
\begin{equation}
w\cdot \caO = s_{\alpha_k}\cdot \left (\cdots s_{\alpha_2} \cdot ( s_{\alpha_1} \cdot \caO)) \right ).
\end{equation}
It is easy to see this is well-defined independent of the choice of reduced expression.  

The following is our main tool to detect failure of rational smoothness, and is a
special case of the results of Springer \cite{S}.

\begin{theorem}
\label{t:springer}
Fix $\caO \in K \bs \caB$ and a closed orbit $\caO_\cl \in K\bs \caO$.  Consider
\[
S(\caO, \caO_\cl) = 
\{\alpha \in \Delta^+\; | \; s_\alpha \cdot \caO_\cl \neq \caO_\cl \text{ and } s_\alpha \cdot \caO_\cl \subset \ol{\caO}\}.
\]
If
\begin{equation}
\label{e:springer}
\#S(\caO,\caO_\cl) > \dim(\caO) - \dim(\caO_\cl),
\end{equation}
then the closure of $\caO$ is not rationally smooth.
\end{theorem}

\begin{remark}
\label{r:springer}
If $\frg$ contains a Cartan subalgebra which is fixed pointwise by $\theta$, then 
the condition $s_\alpha \cdot \caO_\cl \neq \caO_\cl$ is, in the usual terminology, equivalent to $\alpha$ being noncompact imaginary for  $\caO_\cl$.  
See \cite[4.3]{RS}, for instance.
\end{remark}

\medskip

The necessary condition for rational smoothness furnished by Theorem \ref{t:springer}
is a priori rather weak; for example, the analogous necessary condition for rational smoothness of complex Schubert varieties is far from sufficient, even in low rank.  We will see below, however, that this condition is both necessary and sufficient for rational smoothness in all the cases that we consider.
(Even though we make no use of them, we mention that a number of other powerful techniques exist
 for detecting (rational) smoothness; see \cite{brion}, \cite{brion2},  and \cite{CK}, for example, and the references
 therein.)

\section{$\U(p,q)$}
\label{s:upq}
We specialize to the setting of \cite{mcg:upq} before returning to questions
of isogeny at the end of this section.
Fix integers $p\+q = n$, and a signature $(p,q)$
Hermitian form $\langle \cdot, \cdot \rangle$ on an $n$-dimensional
complex vector space $V_\C$.
Let $G = \GL(V_\C)$ (the invertible complex linear endomorphisms of $V_\C$
with determinant one)
and let $G_\R$ denote the subgroup
of $G$ preserving $\langle \cdot, \cdot \rangle$.
Write $\theta$ for the complexification of a Cartan involution of $G_\R$
and set $K = G^\theta$.  Then $K\simeq \GL(p,\C)\times
\GL(q,\C)$.
Choose coordinates so that $\Delta^+ = \{e_i -e_j \; | \; i < j\}$.

In the present setting, the twisted involutions of \cite[Section 3]{RS} 
parametrizing $K \bs \caB$ amount to {\em involutions of
$S_{p+q}$ with signed fixed points of signature $(p,q)$}; that is,
involutions in the symmetric group $S_{p+q}$
whose fixed points are labeled with
signs (either $+$ or $-$) so that half the number of non-fixed points plus
the number of $+$ signs is exactly $p$ (or, equivalently, half the number
of non-fixed points plus the number of $-$ signs is $q$). 
The parametrization is arranged to have the following feature.
Suppose $\caO$ is parametrized by an involution with signed fixed points
whose 
underlying 
involution in the symmetric group is $\sigma$.  Then there is a representative 
$\frb_\caO = \frh_\caO \oplus \frn_\caO$ (with $\frh_\caO$ $\theta$-stable)
with the following property.  Write $\Delta_\caO^+$ for the
roots of $\frh_\caO$ in $\frg$.
There
is a unique inner automorphism of $\frg$ carrying $\frh$ to $\frh_\caO$ and
$\Delta^+$ to $\Delta^+_\caO$.  Using it we may
transport the action of $\theta$ on $\Delta^+_\caO$ to our fixed system
$\Delta^+$.  Once this is done, the key property is that
\begin{equation}
\label{e:thetaupq}
\theta(e_i - e_j) = e_{\sigma(i)} - e_{\sigma(j)}.
\end{equation}

Let $\Spq$ denote the set of signed involutions of $S_{p+q}$  with
signature $(p,q)$ as in the previous paragraph.
For the purposes of formulating pattern avoidance results, we introduce
the following notation.
An element $\gamma \in \Spq$ 
will be identified with an $n$-tuple $(c_1,\dots, c_{n})$, with each $c_i$
either a natural number, a $+$, or a $-$ such that: every natural
number occurs exactly twice; and the number of distinct natural
numbers plus the number of $+$ entries is exactly $p$.  For later use
in Definition \ref{d:po} below, we say that two such strings are
{\em equivalent} if they have the same signs in the same position and pairs
of equal numbers in the same positions.   (So $11+-22$ is equivalent to $22+-33$
and $55+-22$, for instance.)  In any event, the correspondence between
equivalence classes of such strings and involutions with signed fixed points is clear: pairs
of equal natural numbers $c_i = c_j$ in the string correspond to indices $i$ and $j$ interchanged
by the involution, and a sign $c_i$ in the string corresponds to a label of the fixed points $i$
of the involution.
We will generally not distinguish
between elements of $\Spq$ and (equivalence classes of) such strings (which, we remark,
are called ``clans'' in \cite{yamamoto}).

We now turn to an explicit description of the closure order.
(It may
be helpful to refer to Figure 1
when reading
the discussion below.)  
Fix $\gamma = (c_1, \dots,c_n) \in \Spq$,
a simple root $\alpha = e_i - e_{i+1} \in \Delta^+$,
and recall the action (described before
\eqref{e:thetaupq}) of $\theta$ on $\alpha$
determined by $\caO_\gamma$.
Then
\[
\dim\left ( s_\alpha\cdot  \caO_\gamma \right )
= 1 + \dim(\caO_\gamma)
\]
if and only if one of the following conditions hold
(using terminology as in \cite[4.3]{RS}, for instance),
\begin{enumerate}

\item[(a)] $\alpha$ is complex (i.e.~$\theta(\alpha) \neq \pm\alpha$)
and $\theta\alpha \in \Delta^+$;

\item[(b)] $\alpha$ is imaginary (i.e.~$\theta(\alpha) = \alpha$) and
noncompact.
\end{enumerate}
Using \eqref{e:thetaupq} for the former condition, and a calculation in 
$\U(1,1)$ for the latter, these may be formulated as the following
conditions on positions $i$ and
$i+1$ in $\gamma$:
\begin{enumerate}
\item[(a)] $c_i$ and $c_{i+1}$ are unequal natural numbers such
that $j<k$ where $c_i = c_j$ and $c_{i+1} = c_k$ (and $j\neq i$
and $k \neq i+1$);

\item[{($\text{a}'$)}] 
$c_i$ is a sign, $c_{i+1}$ is a number and the entry
$c_k$ with $c_k = c_{i+1}$ ($k \neq i+1$) satisfies $i<k$;

\item[{($\text{a}''$)}]
 $c_i$ is a number, $c_{i+1}$ is a sign, and
the other entry
$c_k$ with $c_k = c_i$ ($k \neq i$) satisfies $k<i+1$;

\item[(b)] $c_i$ and $c_{i+1}$ are opposite signs.
\end{enumerate}
Thus if ${i=2}$, 
then
$\gamma _1= (1,1,2,2)$ satisfies the first condition above, 
$\gamma_2=(-,+,1,1)$
satisfies the second, $\gamma_3 = (1,1,-,+)$ the third,
$\gamma_4 = (1,+,-,1)$ the fourth,
and $(1,2,1,2)$ or $(1,+,1,-)$ satisfies none of them.
In each case the dense $K$ orbit in $s_\alpha \cdot \caO_\gamma$
is parametrized by $\gamma' = (c'_1,\dots,c_n') \in \Spq$ 
which differs from $\gamma$ only in the $i$ and $i+1$ entries in
each respective case
as follows:
\begin{enumerate}
\item[(a)]
$c'_i = c_{i+1}$ and $c'_{i+1}= c_i$.

\item[(b)] $c_i' = c_{i+1}'$ is a natural number.
\end{enumerate}
So for the examples listed above, we have 
$\gamma_1' = (1,2,1,2)$, $\gamma_2' = (-,1,+,1)$,
and $\gamma_3' = (1,-,1,+)$,
and $\gamma'_4 = (1,2,2,1)$.  

The previous paragraph thus gives a complete description of the
weak closure order.  After applying the recursive procedure
given above, one obtains an explicit description of the full 
closure order.  
In particular,
all of the following operations move from one element $\gamma
\in \Spq$ to a higher one in
the order:  replace a pair of (not necessarily adjacent) opposite
signs by a pair of equal numbers; or interchange a number with a sign
(again not necessarily adjacent to it) so as to move the number
farther away from its equal mate in the string (and on the same side);
or interchange a pair $a,b$ of equal numbers with $a$ to the left of
$b$ provided that the mate of $a$ lies to the left of the mate of $b$.
Thus (the orbit corresponding to) $(1,+,1,-)$ lies below $(1,2,1,2)$
and $(1,+,-,1)$, while 
$(1,2,1,3,2,3)$ lies below $(1,3,1,2,2,3)$ but not below $(1,3,1,3,2,2)$.
In particular, the closed orbits are parametrized by elements of $\Spq$ consisting
only of signs, while the open orbit is parametrized by $(1,2,\dots, 2q, +,
\dots , +,  2q-1,2q,2q-3,2q-2,\dots,1,2)$, 
with $2p - 2q$ plus signs, if $p > q$. 

By similar considerations,  one quickly deduces the following dimension formula
(as in \cite[Section 2.3]{yamamoto}, for instance) for the orbit $\caO_\gamma$
parametrized by $\gamma = (c_1,\dots,c_n)$.  Let
\begin{equation}
\label{e:length}
l(\gamma) = \sum_{c_i=c_j \in \bbN, i<j}\left( j-i-
\#\{k \in \bbN \; | \; c_s = c_t = k \text{ for some }
s < i<t<j \} \right ).
\end{equation}
Then
\[
\dim(\caO_\gamma) = 
d(K) +  l(\gamma),
\]
where $d(K)$ is the dimension of the flag variety for $K$, namely
$\frac12(p(p-1) + q(q-1))$.

\begin{defn}
\label{d:po}
We say that an involution with signed fixed points $(c_1, \dots, c_n) \in \Spq$ includes the pattern
$(d_1, \dots, d_m)$ if there are indices $i_1 < \cdots < i_m$
so that the (possibly shorter) string $(c_{i_1}, \dots, c_{i_m})$
is equivalent to $(d_1, \dots, d_m)$.  We say that $\gamma$
avoids $(d_1,\dots, d_m)$ if it does not include it.
For instance, $(1,1,2,+,3,2,-,3)$:  contains the pattern $(1,1,2,2)$ (by considering
$i_1 = 1, i_2=2, i_3 =3,$ and $i_4=6$), contains
$(1,2,1,2)$ (by considering
$i_1 = 3, i_2=5, i_3=6,$ and $i_4=8$), and contains $(1,+,1,-)$
 (by considering
$i_1 = 3, i_2=4, i_3=6,$ and $i_4=7$);
but avoids $(1,+,-,1)$, $(1,+,+,1)$, and $(1,2,2,1)$.
\end{defn}
\medskip

\noindent 
Here is the main result from \cite{mcg:upq}.

\begin{theorem}
\label{t:upq}
Fix $\gamma \in \Spq$, an involution in $S_{p+q}$ with signed fixed points of signature
$(p,q)$ (as defined above).  If $\gamma$
includes one of the patterns
$(1,+,-,1)$, $(1,-,+,1)$, $(1,2,1,2)$,  $(1,+,2,2,1)$, $(1,-,2,2,1)$, $(1,2,2,+,1)$, 
$(1,2,2,-,1)$, or $(1,2,2,3,3,1)$, then the closure of
$\caO_\gamma$ is not rationally smooth.
In all other cases,  $\bar\caO_\gamma$ is the support of a derived functor module and so has smooth closure.
In particular, an orbit has smooth closure if and only if it has rationally
smooth closure, or if and only if the condition of Theorem \ref{t:springer} fails for some closed orbit below it.
\end{theorem}

\bigskip
Next we turn to issues of isogeny.  First, we switch notation and consider the simply connected simple group $G= \SL(n,\C)$.
Since the center of $\GL(n,\C)$ acts trivially on $\caB$, the above discussion applies without change for symmetric subgroup $K = \mathrm{S}(\GL(p,\C) \times \GL(q,\C))$ of $G$.  Now let $\ol G$ be a quotient of $G$ by a subgroup $F < \bbZ/n$ of the
center of $G$.  Let $\ol K$ be the symmetric subgroup of $\ol G$ with Lie algebra $\frs(\frg\frl(p,\C) \oplus \frg\frl(q,\C))$.  Then the
orbits of $\ol K$ on $\caB$ coincide with those of $K$ {\em except} in one case: if $n=2m$, $p=q=m$, and $F$ contains the index $m$ subgroup
of $\bbZ/2m$.  In this case, $\ol K$ has two connected components, and indeed $\ol K$ orbits can be disconnected.  More
precisely, $\ol K$ orbits are parametrized by $\Spq$ modulo the equivalence relation generated by $\gamma \sim -\gamma$,
where $-\gamma$ is obtained from $\gamma$ by reversing all signs of fixed points.
(So, for instance, $(1,+,-,1)$ is equivalent to $(1,-,+1)$, but $(1,2,1,2)$ is equivalent to only itself.)  Equivalence classes
thus have either one or two elements.  In general, if $\ol \gamma$ denotes such a class, the $\ol K$ orbit $\caO_{\ol \gamma}$
parametrized by $\ol \gamma$ breaks into $K$ orbits as
\begin{equation}
\label{e:breaks}
\caO_{\ol \gamma} = \caO_\gamma \; \bigcup \; \caO_{-\gamma}
\end{equation}
which is disconnected if $\gamma \neq -\gamma$.  So the subtleties alluded to in the introduction could
potentially come into play.  In fact, they do not, and one way to see this is to examine the proof of Theorem 
\ref{t:upq}.  It consists of two steps: first proving that if $\gamma$ contains one of the eight indicated patterns,
then it is not rationally smooth; and, second, proving that the remaining orbits are smooth.  The first
is done by finding a suitable closed orbit $\caO_\cl$ in the closure of $\caO_\gamma$ so that Theorem
\ref{t:springer} applies.  It turns out that the identical argument can be carried out for $\ol K$ orbits; that is, 
roughly speaking, Springer's criterion is insensitive to isogeny in this case.   (This need {\em not} always
be true --- see Example \ref{e:bad}.)  The second step can also be carried out
in an analogous manner which, once again, is a special feature of this case.  In the end, one concludes
{\em Theorem \ref{t:upq} holds for the orbits of any symmetric subgroup $\ol K$ of $\ol G$ (isogenous to $\SL(n,\C)$) with Lie
algebra $\frs(\frg\frl(p,\C) \oplus \frg\frl(q,\C))$}.  

\bigskip

We conclude this section by recalling (as in \cite[Section 2]{RS}) the natural 
action of $W = W(\frh,\frg)$ on the set of twisted involutions parametrizing
$K$ orbits on $\caB$.  (We will need this for formulating some results
in Section \ref{s:sostar}.)  
On the level of $\Spq$ we write
$w \times \gamma$ for the action of $w$ on $\gamma$.  
Explicitly, it amounts to 
the obvious action of the symmetric group: the
$w(i)$th component of (the string parametrizing) $w\times\gamma$ is simply the $i$th component of
(the string parametrizing) $\gamma$. 

\section{$\Sp(p,q)$}
\label{s:sppq}

We now turn our attention to $G_\R = \Sp(p,q)$, a real form of $G =
\Sp(2n,\C)$, again deferring isogeny questions to the end of this section.
We start with a brief outline of our strategy.  As the definition
recalled below makes evident, $G_\R$ is a subgroup
of $G_\R' = \U(2p,2q)$.   In this section, we will try (as much as possible)
to reduce the study of $K$ orbits on $\caB$,
the flag variety for $G$, to corresponding results in the previous section
for $K' \simeq \GL(2p,\C) \times \GL(2q,\C)$ orbits on $\caB'$, the
flag variety
for $\GL(2p+2q,\C)$.
More precisely, $\caB$ naturally includes into $\caB'$ and (for appropriate
choices of the Cartan involutions in question) a $K'$ orbit on $\caB'$ either
meets $\caB$ in a single $K$ orbit, or else not at all.  Thus the orbits
of $K$ on $\caB$ are parametrized by a subset of the involutions with
signed fixed points $\Sigma_\pm(2p,2q)$ introduced above.   In fact, we quickly check below
that the closure order on $K$ orbits on $\caB$ is simply the 
appropriate restriction of the closure order of $K'$ orbits on
$\caB'$.  Thus the simplest
result that one could hope for is this: the closure of an orbit for $K$
on $\caB$ has (rationally) smooth closure if and only if the corresponding
$K'$ orbit $\caB'$ which it meets has (rationally) smooth 
closure (and recall these latter orbits have been classified in
Theorem \ref{t:upq}).  This is in fact the content of Theorem
\ref{t:sppq} below, apart 
from an easily stated exception treated in Lemma \ref{l:sppq}.

It is worth remarking that since (up to issues of isogeny) any classical real group outside of type A
is a subgroup of an appropriate $\U(p,q)$, one could attempt to mimic
the strategy of the previous paragraph for any such group.  In the next
section, we do so for $\SO^*(2n)$ and obtain similar results.  But for 
the other classical groups, complications arise; 
see Example \ref{e:bad}.

We return to the details of the case of $\Sp(p,q)$.  Let $\bbH$ denote
the quaternions equipped with the standard bar operation
$\overline{a + bi + cj + dk} = a - bi - cj - dk$.  Embed $\bbC$ in
$\bbH$, as usual,
as elements of the form $a+bi$ and write the corresponding isomorphism
$\bbH \simeq \bbC^2$
as $z = \iIm(z) +\jIm(z)j$.
Let $\langle \cdot, \cdot \rangle$ 
denote a signature $(p,q)$ sesquilinear form on an $n$-dimensional
quaternionic (left) vector space $V_\bbH$. 
Define 
\[
\langle \cdot, \cdot \rangle'  = \iIm \circ
\langle \cdot, \cdot  \rangle
\]
and
\[
\langle \cdot, \cdot \rangle''  = \jIm \circ
\langle \cdot, \cdot  \rangle.
\]
Then $\langle \cdot, \cdot \rangle'$ is a nondegenerate Hermitian form
of signature $(2p,2q)$ on the underlying $2n$-dimensional complex vector
space $V_\C$ and
$\langle \cdot, \cdot \rangle''$ is a nondegenerate alternating
form on $V_\C$.

Let $G_\R = \Sp(p,q)$ denote the subgroup of $\GL(V_\bbH)$  (the group
of invertible
left $\bbH$-linear endomorphisms of $V_\bbH$) preserving
$\langle \cdot, \cdot \rangle$ and view, as we may, $G_\R$ as subgroup 
of $\GL(V_\C)$.  Let $G'_\R$ denote the subgroup
of $G':=\GL(V_\C)$ preserving
$\langle \cdot, \cdot \rangle'$, and finally
let $G$ denote the
subgroup of $\GL(V_\C)$ preserving 
$\langle \cdot, \cdot \rangle''$.  Then $G_\R$ is a real form
of $G \simeq \Sp(2n,\C)$, $G_\R' \simeq \U(2p,2q)$, and
\[
G_\R = G_\R' \cap G.
\]
We adopt the notation of Section \ref{s:upq} for $G_\R'$, 
adding a prime everywhere as appropriate.
Any Cartan involution for $G_\R'$ restricts to one for
$G_\R$.  If we write $\theta'$ for the corresponding complexified
involution of $G'$ and $\theta$ for its restriction to $G$, we naturally
have
\[
\Sp(2p,\C) \times \Sp(2q,\C) \simeq K := G^\theta < K' := (G')^{\theta'} \simeq
\GL(2p,\C) \times \GL(2q,\C).
\]
Recall the natural inclusion of the flag variety $\caB$ for $G$ into $\caB'$,
the flag variety for $G'$.
By the remarks above, in order to classify $K$ orbits on $\caB$
it suffices to determine which orbits $\caO'$ meet $\caB$ nontrivially.
If $\gamma' = (c_1', \dots, c_{2n}')  \in \Sigma_\pm(2p,2q)$
(with notation as in the last section), we say that it is symmetric if
\begin{enumerate}
\item[(i)]
If the entry $c'_i$ is a sign, then $c'_{2n+1-i}$ is the same sign.

\item[(ii)] If $c'_i=c'_j$ are natural numbers, 
then $j\neq 2n+1-i$ and $c'_{2n+1-i} = c'_{2n+1-j}$.
\end{enumerate}
We write  $\SPpq$ for the subset of symmetric elements
in $\Sigma_\pm(2p,2q)$.
Then the $K'$ orbit $\caO'_{\gamma'}$
meets $\caB$ if and only if $\gamma'$ is symmetric 
(\cite{MO}, \cite[Section 4.3]{yamamoto}), and thus
$\SPpq$ parametrizes $K$ orbits on $\caB$.  
Our
next task is to describe the closure order explicitly.

Fix, as in Section \ref{s:upq}, a $\theta'$-stable Cartan subalgebra
$\frh'$, a choice of positive roots, and (in appropriate coordinates)
write $(\Delta')^{+} = \{e'_i - e'_j \; | \; i<j\}$.  Then $\frh:=
\frg \cap \frh'$
is a $\theta$-stable Cartan subalgebra of $\frg$.
Restriction defines a positive system of roots of $\frh$ in $\frg$,
$\Delta = \{2e_i \; | \; 1 \leq i \leq n\} \cup \{e_i \pm e_j \; 1
\leq i < j \leq n\}$.   
Let $B$ denote
the corresponding Borel subgroup of $G$.  Fix a simple root $\alpha$
in $\Delta^+$ and let $S'(\alpha)$
the set of roots in $(\Delta')^{+}$ which restrict to $\alpha$.
Concretely, if $\alpha = e_i -e_{i+i}$,
then $S'(\alpha) = \{\alpha',\alpha''\}$ with $\alpha' = e'_i
\-e'_{i+1}$ and $\alpha'' = e'_{2n-i} -e'_{2n-i+1}$; and if $\alpha =
2e_1$, then $S'(\alpha)= \{\alpha'\}$ where $\alpha'=e_n'-e'_{n+1}$.

Given a simple root $\alpha \in \Delta^+$ 
and $\caO$ in $K \bs \caB$, define $s_\alpha\cdot \caO$ as in \eqref{e:pnot}.
For $\beta \in (\Delta')^+$ and $\caO'$ in $K' \bs \caO'$, define
$s'_\beta\cdot \caO'$ similarly.
The parametrization of
$K$ and $K'$ orbits satisfies the following key geometric
compatibility condition
(which follows easily from unraveling the definitions).
Given any symmetric
$\gamma$, let $\caO_\gamma$ denote the corresponding $K$ orbit on $\caB$
and $\caO'_\gamma$ the corresponding $K'$ orbits on $\caB'$.  Then
\begin{equation}
\label{e:sp1}
\dim\left (  s_\alpha\cdot \caO_\gamma \right) 
= 1 + \dim(\caO_\gamma)
\end{equation}
if and only if
\begin{equation}
\label{e:sp2}
\dim\left (s'_\beta\cdot \caO'_\gamma \right) 
= 1 + \dim(\caO'_\gamma)
\end{equation}
for some (equivalently, any) root $\beta \in S'(\alpha)$.
Moreover, in this case, 
the dense $K'$ orbit in
\begin{equation}
\label{e:product}
\left [ \prod_{\beta \in S'(\alpha)} s'_\beta \right] \cdot 
\left (\caO'_{\gamma}\right ).
\end{equation}
intersected with $\caB$ is the dense $K$ orbit in 
\begin{equation}
\label{e:product2}
\pi_\alpha^{-1}\left ( \pi_\alpha (\caO_\gamma) \right ).
\end{equation}
Note that if the product in \eqref{e:product} has more than one term,
the order of the terms does not
matter since, in this case, the roots in $S'(\alpha)$ are orthogonal.

Using these facts, it is a simple matter to write down the weak closure
order on the level of $\SPpq$ from the description of the
weak closure order described in the previous section. (This is done
in \cite[4.4]{yamamoto}, for instance.)  By examining the
recursive procedure to generate the full closure order from the weak
order, one concludes
that the closure order of $K$ orbits on the level 
$\SPpq$ is
simply the restriction of the order 
on $\Sigma_\pm(2p,2q)$ given in the
previous section.   The dimension of the orbit
parametrized by $\gamma \in \SPpq$ is also easy to read off, and is given by
\[
\dim(\caO_\gamma) = d(K) +
(1/2)\left ( l(\gamma) + \{t \in \bbN \; | \;
 c_s = c_t \in \bbN  \text{ with } s \leq n < t \leq 2n + 1 - s\}\right )
\]
where $d(K)$ is the dimension of the flag variety for $K$, namely $p^2
+ q^2$, and
$l(\gamma)$ is as in \eqref{e:length}.  In particular, closed orbits
are once again
parametrized by elements consisting only of signs, while the open orbit
is parametrized by
\begin{equation}
\label{e:open}
\gamma_\circ(p,q) := (1, 2, \dots, 2q,{+, \dots, +}, 2q-1, 2q,
2q-3,2q-2,\dots,1,2)
\end{equation}
with $p-q$ plus signs
if $p \geq q$, and similarly if $q \geq p$.  
See Figure 2
for a detailed example.

We need some notation for
the next result.  If $\gamma_1 = (c_1, \dots, c_k)$ represents
an involution with signed fixed points,
then let $\gamma_1^r$ be the string with the coordinates of $\gamma$ in reverse
order but with each pair of equal numbers changed to a different pair of numbers.  For example, if
$\gamma_1=(1,2,+,1,2,-)$, then $\gamma_1^r=(-,3,4,+,3,4)$.  For any such
$\gamma_1$, the concatenation $\gamma:=(\gamma_1,\gamma_1^r)$ is
symmetric, so may be viewed as an element of $\SPpq$.

\begin{lemma}
\label{l:sppq}
Fix $\gamma \in \SPpq$.  Suppose that
there are integers $p' + r =p $ and $q' + s = q$ so that $\gamma$ can
be written as the concatenation
\[
(\gamma_1, \gamma_\circ(p',q'),  \gamma_1^r)
\]
with $\gamma_1 \in \Sigma_\pm(r,s)$ such that $\gamma_1$ avoids
the bad patterns of Theorem \ref{t:upq}.
Then the closure of
$\caO_\gamma$ is the support of a derived functor module for $G_\R$ and hence is smooth.  The fiber is isomorphic
to the product of the flag variety for $\Sp(2p'\+2q',\C)$ and the 
closure of the $\GL(r,\C) \times \GL(s,\C)$  orbit parametrized
by $\gamma_1$ (as in Section \ref{s:upq}).
\end{lemma}

\noindent{\bf Proof.}
We start with a general observation.  Suppose $\gamma$ is any
element of $\Spq$
which can be written as a concatenation $(\gamma_1,\gamma_\circ,
\gamma_1^r)$ where
$\gamma_\circ \in \Sigma_\pm(p',q')$
and $\gamma_1 \in \Sigma^\sym_\pm(2r,2s)$.
Set $n' = p'\+ q'$, $n'' =
r\+s$.  Then there
is a $\theta$-stable parabolic subgroup $Q=LU$, unique up to $K$
conjugacy, containing $B$ with 
$L \simeq \Sp(n',\C) \times \GL(n'',\C)$ and 
\[
K \cap L \simeq 
\Sp(p',\C) \times \Sp(q',\C) \times 
\GL(r,\C) \times \GL(s,\C).
\]
The assumption that 
$\gamma = (\gamma_1,\gamma_\circ, \gamma_1^r)$ implies that the
image of $\caO_\gamma$ under projection
$\pi_\caP$ from $\caB \simeq G/B$ to $\caP \simeq G/Q$ is the closed
$K$ orbit of the identity coset $eQ$ and that the fiber of the 
restriction of $\pi_\caP$ to $\caO_\gamma$ is a single orbit of 
$L\cap K$ on $\caB_L$, the flag
variety for $L$.   The fiber is
isomorphic to the product of the $\Sp(p',\C) \times \Sp(q',\C)$ orbit
parametrized
by $\gamma_\circ$ with the orbit of
$\GL(r,\C) \times \GL(s,\C)$ parametrized by $\gamma_1$.  Meanwhile
the $K$ orbit of $eQ$ is closed, and hence isomorphic to a partial
flag variety for $K$.
The closure of $\caO_\gamma$ is thus a fiber bundle over a partial flag variety
for $K$ whose fiber is the closure of the product of the orbits
parametrized by $\gamma_\circ$
and $\gamma_1$.
If we impose
the additional hypothesis that $\gamma_\circ = \gamma_\circ(p',q')$
and $\gamma_1$ avoids the patterns given in
the lemma, then we conclude from Theorem \ref{t:upq} that the orbit closure takes the claimed form.\qed

\begin{theorem}
\label{t:sppq}
Fix $\gamma \in \SPpq$, a symmetric involution of $S_{2p+2q}$ with signed
fixed point of signature $(2p,2q)$ (as defined above).  Suppose $\gamma$ is not of the
form treated
by Lemma \ref{l:sppq} and, further, that $\gamma$ includes one of the bad patterns of Theorem 2.2.  Then the closure of $\caO_\gamma$ does not have rationally smooth closure.  

In all other cases,
the closure of $\caO_\gamma$ is the support of a derived functor module, so is smooth.
In particular, $\caO_\gamma$ has rationally smooth closure if and only if
it has smooth closure, or if and only if the condition of Theorem \ref{t:springer} fails 
for some closed orbit below it.
\end{theorem}  

\noindent{\bf Proof.} If
$\gamma$ is of the form treated by Lemma \ref{l:sppq}, there is
nothing to prove.
So suppose this is not the case and that 
$\gamma$ avoids the patterns given in the theorem.  Then
$\gamma$ takes the form $(\gamma_1, \gamma_1^r)$ with $\gamma_1
\in \Spq$ which avoids the same patterns. 
Thus Lemma \ref{l:sppq} applies (with $p'=q'=0$) to give the required
assertion.

Suppose $\gamma$ contains one of the patterns
of the theorem.  Then there are various geometric ways to 
deduce the failure of the closure of $\caO_\gamma$ to be rationally
smooth by embedding the corresponding failure for the smaller rank group
where, roughly speaking, the pattern resides.  (The exceptions of Lemma \ref{l:sppq}
shows that some care is required.)
A concise (and convenient) way to organize the geometric reduction
is to use Theorem \ref{t:springer}, and this
is how we shall proceed.

In the present context
a root $e_i - e_j$ with $i<j$ is noncompact imaginary for an
orbit $\caO_\cl$ parametrized $\cl$ if
the $i$th and $j$th coordinates of $\cl$ are opposite signs and the root
$e_i + e_j$ is noncompact imaginary if the $i$th and $2n +1
- j$th coordinates are 
opposite signs.  (The root $2e_i$ is never noncompact imaginary.) 
Then 
if $\alpha=e_i-e_j$,
$s_\alpha\cdot\caO_\cl$ is parametrized 
by the element obtained from $\cl$ by
replacing its $i$th, $j$th and $2n+1-j$th, $2n+1-i$th 
coordinates by different pairs of equal numbers;
if instead $\alpha=e_i+e_j$,
the new element is obtained by
replacing the
$i$th, $2n+1-j$th and $j$th, $2n+1 -i$th coordinates by two different pairs of
equal numbers.
As we
remarked above, closed orbits are parametrized by elements
consisting of only signs.  Since we have given the closure order and
dimension formula above, the  criterion of \eqref{e:springer}
becomes very explicit, and we may apply it directly as follows.

Fix a $\gamma=(c_1,\ldots,c_{2n}) \in \SPpq$ including one of the bad
patterns and not of the form treated by Lemma \ref{l:sppq}.  We first
produce a suitable closed orbit $\caO_\cl$ lying below $\caO_\gamma$.
Look first at the natural numbers occurring twice among
$c_1\ldots,c_n$.  Replace the first occurrence of all such numbers by
$+$ and the second by $-$.  Then look at all the natural numbers
occurring just once among $c_1,\ldots,c_n$. Whenever $c_i=c_{2n+1-j}$
is a natural number and $i<j<n$ (so that $c_j = c_{2n+1-i}$ is also a
natural number), replace $c_i$ by $+$ and $c_j$ by $-$.  Finally,
replace $c_{n+1},\ldots,c_{2n}$ by signs in such a way that the
result $\cl=(\cl_1,\ldots,\cl_{2n})$ is symmetric.

Assume for the moment that there are no indices $i<j\le n$ with $c_i$
and $c_j$ equal natural numbers.  Enumerate the indices $i\le n$ with
$c_i$ and $c_{2n+1-j}$ equal natural numbers for some $j\le n$ as
$i_1<\cdots<i_{2m}$.  Define a new element
$\gamma'=(c_1',\ldots,c_{2n}')$ by decreeing that
$(c'_{i_1},\ldots,c'_{i_1+2m-1}) = (1,\ldots,2m)$, that
$(c'_{2n+1-(2m-1)-i_1},\ldots,c'_{2n+1-i_1}) =
(2m-1,2m,2m-3,2m-2,\ldots,1,2)$, and finally that there be as many $+$
signs among $(c'_1,\ldots,c'_n)$ as among $(c_1,\ldots,c_n)$ and
similarly for $-$ signs.  Then $\gamma'$ parametrizes an orbit
$\caO_{\gamma'}$ whose dimension is at least that of $\caO_\gamma$,
but for which the left-hand side of \eqref{e:springer} for
$\caO_{\gamma'}$ (and $\caO_\cl$)
is no larger than it is for $\caO_\gamma$ (and $\caO_\cl$).  Thus if
Springer's criterion implies that the closure of $\caO_{\gamma'}$
fails to be rationally smooth,
or $\caO_{\gamma'}$ has dimension larger than $\caO_\gamma$,
the same is true for the closure of
$\caO_\gamma$.  Note that both sides of \eqref{e:springer} are unaffected if
we replace $\gamma'$ by the new (smaller) element obtained by deleting
all the initial signs $c'_1, \dots, c'_{i_1-1}$ and terminal signs
$c'_{2n-i_1+2}, \dots, c'_{2n}$ from $\gamma'$, and 
replace $\cl$ by the element obtained by deleting the corresponding
initial and terminal entries from $\cl$.  So we may assume $i_1 = 1$.  

We argue inductively based on the number of distinct natural
numbers $k$ appearing in $\gamma'$ as follows.  We will reduce
our analysis to the case of $k=2$, so assume first that $k>2$.  
Consider the new element $\gamma''$ (of smaller size) obtained
from $\gamma'$ first by deleting the initial entries $c'_1, \dots, c_{i_2}'$
and terminal entries $c'_{2n+1-i_2}, \dots, c'_{2n}$,
and then deleting all the initial and terminal signs from the result.
Let $\cl''$ denote the element obtained by deleting the corresponding
entries from $\cl$.
Then in passing from $\caO_{\gamma'}$ to $\caO_{\gamma''}$
(and from $\caO_\cl$ to $\caO_{\cl''}$), 
the difference of the left- and right-hand sides of
\eqref{e:springer} will strictly increase in all cases unless $i_2 =
2$, in which case the difference will remain the same.  Continuing
this procedure, we are thus
led to investigating the ``base case'' of the form $\gamma'' =
(1,2,\epsilon_3,
\dots, \epsilon_n, \epsilon_n, \dots, \epsilon_3, 1, 2)$
and $\cl'' = (+, -, \epsilon_3, \dots, \epsilon_n, 
\epsilon_n, \dots, \epsilon_3, -, +)$.  
There are
two possibilities.  If the signs $\epsilon_i$ are not all the same, one checks
directly that \eqref{e:springer} holds for the base case.  Since each
of the steps leading to the base case weakly increases the difference
of the left- and right-hand sides of \eqref{e:springer}, one concludes that
\eqref{e:springer} holds for $\caO_{\gamma'}$ and $\caO_\cl$, hence for
$\caO_\gamma$ and $\caO_\cl$,
and hence the closure of $\caO_\gamma$ is not rationally smooth, as claimed.  
The other possibility is that all of the signs in the base case are the
same.  In this case, the two sides of \eqref{e:springer} corresponding
to the base case are actually
equal.
Again since each step in the reduction
weakly increases the difference
of the left- and right-hand sides of \eqref{e:springer}, the inequality in \eqref{e:springer} holds for $\caO_\gamma'$ and
$\caO_\cl$ (and hence the closure of
$\caO_\gamma$ fails to be rationally smooth)
except possibly in just one
case, namely when each step reducing to the base 
case does not increase
the difference of the left- and right-hand sides, and when the base case
turns out to have all signs the same.  But the only way this can
happen is if $\gamma'$ is of the form $\gamma_\circ(p',q')$ (possibly
flanked by a number of signs).
If this is the case, the construction of $\gamma'$ shows that
$\gamma$ must be of the same form up to a
permuation of the natural number entries.  Since we have assumed $\gamma$ is
not of the form treated by Lemma 3.1, one checks that the initial passage
from $\gamma$ to $\gamma'$ does indeed increase the difference of the left-
and right-hand sides of (5).  So once again the closure of
$\mathcal{O}_\gamma$ is not rationally smooth, as desired.

 Finally return to the case where 
there are indices $i<j\le n$ with $c_i$ and $c_j$ equal natural
  numbers in $\gamma$, and replace $c_i,c_j$ for every such pair $i<j$ by signs
  as in the definition of $\cl$ above (and define
  $c_{2n+1-i},c_{2n+1-j}$ so that the resulting element is symmetric).  We
  obtain an element parametrizing
  an orbit that either does not contain any of the bad patterns, or
  whose closure fails to be rationally smooth by the above argument.
  Changing the $i$th and $j$th coordinates of this element back to
  $c_i,c_j$ (and similarly for the $2n+1-i$th and $2n+1-j$th
  coordinates), we find that the dimension of the corresponding orbit
  increases, but so too does the left-hand side of
  \eqref{e:springer} by at least the same number, and by a larger
  number if two of the signs between $c_i$ and $c_j$ are different.
We conclude in all cases that the
closure of $\caO_\gamma$ is not rationally smooth, as desired.\qed

\medskip
Finally we must consider the situation for the symmetric subgroup $\ol K$ of
$\ol G = \mathrm{PSp}(2n,\C)$ with Lie algebra $\frs\frp(2p,\C) \times \frs\frp(2q,\C)$.
Then $\ol K$ is connected (and its orbits on $\caB$ coincide with those of $K$)
{\em unless} $p =q$, in which case $\ol K$ has two components.  In this case,
orbits of $\ol K$ are parametrized by equivalence classes in $\SPpq$ generated
for the relation generated by $\gamma \sim -\gamma$.  Just as in the case discussed
at the end of Section \ref{s:upq},
 if $\ol \gamma$ denotes such a class, the $\ol K$ orbit $\caO_{\ol \gamma}$
parametrized by $\ol \gamma$ breaks into $K$ orbits as
\[
\caO_{\ol \gamma} = \caO_\gamma \; \bigcup \; \caO_{-\gamma}
\]
which is disconnected if $\gamma \neq -\gamma$.  Once again one may retrace
the steps of the proof of Theorem \ref{t:sppq} to show that no complications
arise.  One thus concludes: 
{\em Theorem \ref{t:sppq} holds for the orbits of any symmetric subgroup $\ol K$ of $\ol G$ 
(isogenous to $\Sp(2p+2q,\C)$) with Lie
algebra $\frs\frp(2p,\C) \oplus \frs\frp(2q,\C)$}.

\section{$\SO^*(2n)$}
\label{s:sostar}
We follow the same strategy as outlined at the beginning of Section
\ref{s:sppq}.
Let $\langle \cdot, \cdot \rangle$ denote a nondegenerate skew sesquilinear
form on an $n$-dimensional quaternionic vector space $V$, so that
\[
\langle u, v \rangle = -\overline{\langle v, u \rangle}.
\]
As in the
previous section, set
\[
\langle \cdot, \cdot \rangle'  
= \iIm \circ
\langle \cdot, \cdot  \rangle
\]
and
\[
\langle \cdot, \cdot \rangle''  = \jIm \circ
\langle \cdot, \cdot  \rangle.
\]
Then $\langle \cdot, \cdot \rangle'$ is a nondegenerate Hermitian form
of signature $(n,n)$ on the underlying $2n$-dimensional complex vector
space $V_\C$ and
$\langle \cdot, \cdot \rangle''$ is a nondegenerate symmetric
form on $V_\C$.
Let $G_\R = \SO^*(2n)$ denote the subgroup of $\GL(V_\bbH)$ preserving
$\langle \cdot, \cdot \rangle$ and view $G_\R$ as subgroup 
of $\GL(V_\C)$.  Let $G'_\R$ denote the subgroup
of $G':=\GL(V_\C)$ preserving
$\langle \cdot, \cdot \rangle'$, and
let $G_\pm$ denote the
subgroup of $\GL(V_\C)$ preserving 
$\langle \cdot, \cdot \rangle''$.  Then
$G_\pm \simeq \O(2n,\C)$, $G_\R' \simeq \U(n,n)$, and
\[
G_\R = G_\R' \cap G_\pm.
\] 
(In fact, every element of $G_\R$ has determinant one, so $G_\R$ is indeed
a real form of the connected algebraic group $G\simeq \SO(2n,\C)$ consisting
of determinant one elements in $G_\pm$.) 
Fixing compatible involutions $\theta$ and $\theta'$ as in Section
\ref{s:sppq}, 
we naturally
have
\[
\GL(n,\bbC) \simeq K := G^\theta < K' := (G')^{\theta'} \simeq
\GL(n,\C) \times \GL(n,\C).
\]
Let $\caB$ denote the flag variety for $G$ and $\caB'$ the flag
variety for $G'$.  Once again, there is a natural inclusion
$\caB' \subset \caB$, and if $\caO'$ is an orbit of $K'$ on $\caB'$,
then its intersection with $\caB$ is either empty or consists of a single
orbit of $K$ on $\caB$.  
If $\gamma' = (c_1', \dots, c_{2n}') \in \Sigma_\pm(n,n)$, we
say that it is antisymmetric if:
\begin{enumerate}
\item[(i)]
If the entry $c'_i$ is a sign, then $c'_{2n+1-i}$ is the opposite sign.

\item[(ii)] If $c'_i=c'_j$ are natural numbers, 
then $j \neq 2n+1-i$ and $c'_{2n+1-i} = c'_{2n+1-j}$.

\item[(iii)] Among the entries $c_1',\ldots,c_n'$, the number of $+$
  signs plus the number of pairs
of equal natural numbers is even. 
\end{enumerate}
Write $\SOpq$ for the set of such antisymmetric elements.  
Notice there is a choice of sign in (iii) above.  In Section
\ref{s:upq}, when $p=q=n$,
the parametrization of orbits can obviously be twisted by an outer automorphism
of $\U(n,n)$, the effect of which is to change
all signs. 
Possibly after twisting the parametrization of Section \ref{s:upq}, 
the $K'$ orbit $\caO'_{\gamma'}$
meets $\caB$ if and only if $\gamma'$ is antisymmetric.
Thus $\SOpq$ parametrizes $K$ orbits on $\caB$ \cite{MO}.
(It may be helpful to refer to Figures 3 and 4
in the course of the discussion below.  In Figure 3 we twist condition (iii) above, requiring that the number of $-$ signs and pairs of equal natural numbers among the first three entries in our clans is even.)

Fix a $\theta'$ stable Cartan subalgebra
$\frh'$, a choice of positive roots, and
write $(\Delta')^{+} = \{e'_i - e'_j \; | \; i<j\}$.  Then $\frh:=
\frg \cap \frh'$
is a $\theta$-stable Cartan subalgebra of $\frg$.
Choose a positive subset of roots of $\frh$ in $\frg$, in appropriate
coordinates write $\Delta^+ = \{(e_i \pm e_j) \; | \; i < j\}$, and
let $B$ denote
the corresponding Borel subgroup of $G$.  Fix a simple root $\alpha$
in $\Delta^+$ and let $S'(\alpha)$
the set of roots in $(\Delta')^{+}$ which restrict 
to $\alpha$.  
For instance, if $\alpha = e_i -e_{i+i}$,
then $S'(\alpha) = \{\alpha',\alpha''\}$ with $\alpha' = e'_i
\-e'_{i+1}$ and $\alpha'' = e'_{2n-i} -e'_{2n-i+1}$.

Adopt notation analogous to that around Equation \eqref{e:sp1},
fix $\alpha = e_i-e_{i+1} \in \Delta^+$, and an antisymmetric element $\gamma$
of signature $(n,n)$.  Then once again we have the identical
conclusions of Equations \eqref{e:sp1}--\eqref{e:product2}.
The situation for $\alpha = e_{n-1} +e_n$ is more subtle, however.
Set $\alpha_{n-1} = e_{n-1} -e_n$, fix an 
antisymmetric element $\gamma$ of signature $(n,n)$, and
let $s'$ denote the reflection in the simple root $e'_n - e'_{n+1}$
in $(\Delta')^+$.  Recall the action described at the end
of Section \ref{s:upq}.
Then 
\begin{equation}
\label{e:sostar1}
\dim\left (s_\alpha\cdot \caO_\gamma\right) 
= 1 + \dim(\caO_\gamma)
\end{equation}
if and only if
\begin{equation}
\label{e:sostar2}
\dim\left (s'_\beta \cdot \caO'_{s' \times\gamma}\right) 
= 1 + \dim(\caO'_{s'\times \gamma})
\end{equation}
for some (equivalently, any) root $\beta \in S'(\alpha_{n-1})$.
Moreover, in this case, 
if $\caO'_\delta$ is the dense $K'$ orbit in
$s'_\beta \cdot \caO'_{\gamma},$
then $\caO'_{s' \times \delta}$ 
intersected with $\caB$ is the dense $K$ orbit in 
$s_\alpha\cdot \caO_\gamma$.
Using this, one may deduce 
that the closure order of $K$ orbits on the level of $\SOpq$
coincides with
the restriction of the order on $\Sigma_\pm(n,n)$ given in
Section \ref{s:upq}.   The dimension of the orbit parametrized by an
$\gamma \in \SOpq$ is given by
\[
\dim(\caO_\gamma) = d(K) + 
(1/2)\left ( l(\gamma) - \{t \in \bbN \; | \;
 c_s = c_t \in \bbN  \text{ with } s \leq n < t \leq 2n + 1 - s\}\right )
\]
where $d(K)$ is the dimension of the flag variety for $K$, namely
$\frac12 n(n-1)$, and
$l(\gamma)$ is as in \eqref{e:length}.  Closed orbits are once again
parametrized by elements consisting only of signs, while the open orbit
is parametrized by 
\begin{equation}
\label{e:open3}
\gamma_\circ(n,n) := 
(1, 2, \dots, 2m-1,2m,2m-1,2m, \dots, 1, 2),
\end{equation}
if $n=2m$ is even and 
\begin{equation}
\label{e:open4}
\gamma_\circ(n,n) := 
(1, 2, \dots, 2m,-,+,2m-1,2m,\dots,1,2),
\end{equation}
if $n=2m+1$ is odd.  We let $\pm \gamma_\circ(n,n)$ denote either
$\gamma_\circ(n.n)$ or the element obtained from inverting the
signs in $\gamma_\circ(n,n)$ (which differs from $\gamma_\circ(n.n)$
only if $n$ is odd and, in this case, no longer satisfies condition (iii)
above).

We record two more operations on $\SOpq$.
Fix $\gamma \in \SOpq$ and
let $s'$ be as above.  Let
$\gamma'$ denote the element obtained  by changing
all signs in $s' \times \gamma$ (described at the end
of Section \ref{s:upq}).  Then exactly one element of
 $\{\gamma',s' \times \gamma'\}$
is antisymmetric.  Let $\tau(\gamma)$ denote this element.
Thus $\tau$ is an involution on $\SOpq$, and hence can be
interpreted as an involution of the set of $K$ orbits on $\caB$.  It
coincides with the action of an outer automorphism 
of $G$.  
(As an example, $\tau$ corresponds to the obvious
symmetry in Figures 3 and 4 below.)
Finally, if $\gamma \in \Sigma_\pm(r,s)$, denote by $\gamma^{-r}$ the element
obtained from $\gamma$
by reversing its coordinates, changing all of its signs, and replacing
every pair of equal natural numbers by a different pair of equal
natural numbers.  Then the concatenation $(\gamma,\gamma^{-r})$ is
antisymmetric and so may be viewed as an element of 
$\Sigma^\asym_\pm(r+s,r+s)$

\begin{lemma}
\label{l:sostar}
Fix $\gamma \in \SOpq$
\begin{enumerate}
\item[(a)]
Suppose that there are
integers $r + s + n' = n$ so that $\gamma$ can be written as the concatenation
\[
(\gamma_1, \pm \gamma_\circ(n',n'),  \gamma^{-r}_1)
\]
where $\pm \gamma_\circ(n',n')$ is defined after \eqref{e:open4},
$\gamma_1 \in \Sigma_\pm(r,s)$  which avoids the bad patterns
of Theorem \ref{t:upq}, and
$\gamma_1^{-r}$ is defined as above.
Then the closure of $\caO_\gamma$ is the support of a derived functor module
and hence is
smooth.  
\item[(b)]
Suppose $n=2m$ is even and that 
$\gamma$ can be written as 
\[
(1, \gamma_1, 2, 1, \gamma_1^{-r}, 2),
\]
where $\gamma_1 \in \Sigma_\pm(r-1.s-1)$ , $r+s =m$, 
avoids the bad patterns of Theorem \ref{t:upq} and in addition avoids the patterns $(\pm,3,3), (3,3,\pm)$, and $(2,2,3,3)$, so that $\gamma$ avoids the bad patterns of Theorem 2.2 apart from $(1,2,1,2)$.  Then the closure of
$\caO_\gamma$ is a fiber bundle with smooth fiber over a 
partial flag variety for $K$,
and hence is smooth.  The fiber is isomorphic to
the closure of the $\GL(r,\C) \times \GL(s,\C)$  orbit parametrized
by the element  $(1, \gamma_1, 1) \in \Sigma_\pm(r,s)$.
\end{enumerate}
\end{lemma}

\noindent{\bf Proof.}
Part (a) is proved much the same way as Lemma \ref{l:sppq} and we omit the
details.  
For (b), note that the closure of the orbit parametrized by $\gamma$ is
isomorphic to the one parametrized by the outer automorphism conjugate 
$\tau(\gamma)$ described above.  But if $\gamma$ has the form indicated
in (b), then $\tau(\gamma)$ has the form indicated in (a).  Thus (b)
follows from (a).  \qed

\begin{theorem}
\label{t:sostar}
Fix $\gamma \in \SOpq$, an antisymmetric involution of $S_{2n}$ with
signed fixed points of signature $(n,n)$ (as defined above).  Suppose
$\gamma$ is not of
the form treated
by Lemma \ref{l:sostar} and, further, that $\gamma$ includes one of the bad patterns
of Theorem 2.2.
Then the closure of $\caO_\gamma$ does not have rationally smooth closure.  

In all other cases, 
the closure of $\caO_\gamma$ is a fiber bundle with smooth fiber over a 
partial flag variety for $K$, and hence is smooth.  
In this case, there are integers $r+s+n' = n$ so that the fiber
is isomorphic to the product of the flag variety for $\SO(2n',\C)$ with
the closure of an orbit of $\GL(r) \times \GL(s,\C)$
on the flag variety for $\GL(r\+s,\C)$ with $r\+s =n$.

In particular, $\caO_\gamma$ has rationally smooth closure if and only if
it has smooth closure, or if and only if the condition of Theorem \ref{t:springer} 
fails for every closed orbit below it.
\end{theorem}
\noindent{\bf Proof.}
This is very similar to the proof of Theorem \ref{t:sppq}.  We omit
the details, apart from noting that the clans $(1,+,3,3,2,1,4,4,-,2)$ and $(1,3,3,4,4,2,1,5,5,6,6,2)$, which satisfy the conditions of Lemma 4.1, apart from the further pattern avoidance condition there, correspond to orbits with rationally singular closure, as seen by applying Theorem 1.1 to suitable closed orbits below them. 
\qed

\medskip

Finally we turn to issues of isogeny.  First assume $n$ is odd, so that the center of $\Spin(2n,\C)$ is
$\bbZ/4$.  So the three complex groups
$\ol G$ to consider are $\Spin(2n,\C)$, $\SO(2n,\C)$, and $\mathrm{PSpin}(2n,\C)$.  It turns
out that the symmetric subgroup $\ol K$ with Lie algebra $\frg\frl(n,\C)$ is connected
in each of these cases, and that the orbits of $\ol K$ on $\caB$ are insensitive to isogeny.
So Theorem \ref{t:sostar} applies without change.

The case of $n=2m$ even is more interesting since the center of $\Spin(2n,\C)$ is $\bbZ/2 \times
\bbZ/2$.  Write $\SO(2n,\C)$ for the quotient by the diagonal $\bbZ/2$, $\SO'(2n,\C) \simeq
\SO''(2n,\C)$ for either of quotients by an off-diagonal $\bbZ/2$, and $\mathrm{PSpin}(2n,\C)$
for the adjoint group. Write $K_\sc$, $K$, $K'$, and $K_\ad$ for the corresponding symmetric
subgroups with Lie algebra $\frg\frl(n,\C)$.  The orbits of $K_\sc$ and $K$  on $\caB$ always
coincide, so are treated by Theorem \ref{t:sostar}.  If $m$ is even, these orbits also coincide with 
the orbits of $K'$ on $\caB$.  If $m$ is odd, the orbits of $K'$ instead coincide with the orbits
of $K_\ad$ on $\caB$.  Finally, the orbits of $K_\ad$ are parametrized by equivalence classes
in $\SOnn$ for the relation generated by $\gamma \sim \tau(\gamma)$ for the involution of
$\SOnn$ described above.  
If $\ol \gamma$ denotes such a class, the $\ol K$ orbit $\caO_{\ol \gamma}$
parametrized by $\ol \gamma$ breaks into $K$ orbits as
\[
\caO_{\ol \gamma} = \caO_\gamma \; \bigcup \; \caO_{\tau(\gamma)}
\]
which is disconnected if $\gamma \neq \tau(\gamma)$.   Nonetheless the (omitted) proof
of Theorem \ref{t:sostar} goes through unchanged to establish that the statement of the theorem
requires no modification in this case.  

The conclusion of the discussion is:
{\em Theorem \ref{t:sostar} holds for the orbits of any symmetric subgroup $\ol K$ of $\ol G$ 
(isogenous to $\Spin(2n,\C)$) with Lie
algebra $\frg\frl(n,\C)$}.

\section{examples}
\label{s:examples}
We conclude with several examples of the closure order of $K$ orbits
on $\caB$.   Vertices in the diagrams below correspond to orbits, and
orbits of the same dimension appear on the same column
(or, in the case of Figure 2, 
row).  
Recall that the closure order is generated by relations in codimension
one; so each diagram need only keep track of such relations.
Dashed edges correspond to relations not present in the weak
closure order, as described at the beginning of Section \ref{s:upq}.

Figure 1 
corresponds to the case of $\SU(2,2)$; that is, the case where $G= \SL(4,\C)$,
$K = \mathrm{S}(\GL(2,\C) \times \GL(2,\C))$, and $\caB$ consists of  complete flags
in $\bbC^4$.
(Without the dashed edges and boxed vertices this is \cite[Figure 7]{MO}.)
The labeling of simple roots is given by
  \[
\xy
(0,3)*{1}+(0,-3)*{\bullet}; (10,3)*{2}+(0,-3)*{\bullet} **\dir{-};
 (10,3)*{2}+(0,-3)*{\bullet}; (20,3)*{3}+(0,-3)*{\bullet} **\dir{-};
\endxy
\]
The boxed vertices correspond to orbits with nonsmooth closures
according to the pattern criterion
given in Theorem \ref{t:upq}.  Suppose now $\ol G$ is a nontrivial
quotient of $G$ by a central subgroup, and let $\ol K$ be the
corresponding symmetric subgroup of $\ol G$.
According to 
the discussion at the end of
Section \ref{s:upq}, quotienting Figure 1 
by the obvious 
$\bbZ/2$ symmetry gives the corresponding picture for $\ol K$ orbits
on $\caB$.

Consider next the case of $\Sp(2,2)$ where $G = \Sp(8,\C)$ and
$K = \Sp(4,\C) \times \Sp(4,\C)$.  
According to the details of Section
\ref{s:sppq}, there are 42 orbits for $\Sp(2,2)$,
too many to fit in a reasonable diagram.
Instead consider $\ol G = \PSp(8,\C) = G/F$ where $F$ is the two
element center of $G$.  Then $\ol K = K/F$ is a symmetric subgroup
of $\ol G$ corresponding to $\PSp(2,2)$.    The orbits of $\ol K$
on the flag variety are given in
Figure 2. 
(This graph, without the 
dashed edges and boxed vertices, is \cite[Figure 15]{MO}.)
Simple roots are labeled as
  \[
\xy
(0,3)*{1}+(0,-3)*{\bullet}; (10,3)*{2}+(0,-3)*{\bullet} **\dir{-};
(10,3)*{2}+(0,-3)*{\bullet}; (20,3)*{3}+(0,-3)*{\bullet} **\dir{-};
(20,3)*{3}+(0,-3)*{\bullet}; (30,3)*{4}+(0,-3)*{\bullet} **\dir{=};
\endxy
\]
According to the discussion at the end of Section \ref{s:sppq},
an orbit labeled by an element of $\Sigma_\pm^\sym(2,2)$ which
contain signs is the (disconnected) union of two $K$ orbits.  
(Nonetheless, such an orbit has smooth closure if and only if each connected
component has smooth closure.)
Boxed vertices in Figure 2 
correspond to $\ol K$ orbits with nonsmooth closures
according to the criterion of Theorem \ref{t:sppq}.

Figure 3 
gives the case of $\SO^*(6)$, i.~e.~$K=\GL(3,\C)$ orbits acting on the flag
variety for $\frs\fro(6,\C) \simeq \frs\frl(4,\C)$.  (This
diagram is \cite[Figure 20]{MO}.)   As noted above, we twist condition (iii), replacing + signs in this condition by - ones. 
Simple roots are labeled as
 \[
\xy
(0,3)*{1}+(0,-3)*{\bullet}; (10,3)*{2}+(0,-3)*{\bullet} **\dir{-};
 (10,3)*{2}+(0,-3)*{\bullet}; (20,3)*{3}+(0,-3)*{\bullet} **\dir{-};
\endxy
\]
Since
$\SO^*(6)$ is a quotient of $\mathrm{SU}(3,1)$ by the subgroup
$\bbZ/2$ (defined over $\bbR$) of the center $\bbZ/4$ of $\mathrm{Spin}(6,\C) = \SL(4,\C)$, 
either Theorem \ref{t:upq} (and the discussion at the end of Section \ref{s:upq}) or Theorem
\ref{t:sostar} applies to give that all orbits have smooth closure in this
case.

Finally, Figure 4 
gives the case of $\SO^*(8)$, i.e.~$\GL(4,\C)$
orbits on the flag variety for $\frs\fro(8,\C)$.
(This diagram without the dashed edges or boxed vertices is
\cite[Figure 19]{MO}.)
Simple roots are labeled as
 
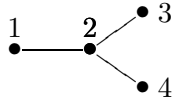
\[
\xy
(0,3)*{1}+(0,-3)*{\bullet}; (10,3)*{2}+(0,-3)*{\bullet} **\dir{-};
 (10,3)*{2}+(0,-3)*{\bullet}; (20,5)*{3}+(-3,0)*{\bullet} **\dir{-};
  (10,3)*{2}+(0,-3)*{\bullet}; (20,-5)*{4}+(-3,0)*{\bullet} **\dir{-};
\endxy
\]
To conserve space, we introduce the following shorthand 
(as in \cite{MO}).  The eight symbols of an element $\gamma \in \Sigma_\pm^\asym(4,4)$
are compressed
to just four.  The signs in the compressed symbol match those
in the first four coordinates of $\gamma$; 
a pair
of numbers in $\gamma$
in positions $i$ and $j$ with $1 \leq i<j \leq 4$ is 
represented by a lower-case  letter in positions
$i$ and $j$ of the compressed symbol; and a pair of numbers
in $\gamma$ in positions $i \leq 4 $ and $9-j > 4$,
is represented by an upper-case letter in positions $i$ and
$j$ of the compressed symbol.  So, for instance, the $\gamma = 1 + - 1 2 + - 1$
becomes $a+-a$; $+12+-12-$ becomes $+AA+$; $12123434$
becomes $abab$; and $12341234$ becomes $ABBA$.
With this convention, the boxed vertices
correspond to orbits with nonsmooth closures
according to the condition of Theorem \ref{t:sostar}.
Isogeny considerations amount (possibly) to folding the figure by the obvious
symmetry, as explained at the end of Section \ref{s:sostar}.

\begin{figure}
\label{f:u22-example}
$${\tiny
{
\xymatrixcolsep{2pc} 
\xymatrixrowsep{3pc}
\xymatrix@d
{
&&&&&1221\\
&&&\boxed{1+-1} \ar@{>}[urr]^2 && \boxed{1212} \ar@{>}[u]^{1,3} && 
\boxed{1-+1}\ar@{>}[ull]^2 \\
&1+1- \ar@{>}[urr]^3  \ar@{.>}[urrrr]_<<<<<<<<{} 
&& +1-1 \ar@{>}[u]_<<<<1 \ar@{.>}[urr]_<<<<<<{} 
&& 1122 \ar@{>}[u]^2 \ar@{.>}[urr]_<<<<<<{}  \ar@{.>}[ull]^<<<<<{}  
&& -1+1 \ar@{>}[u]^<<<<<1 \ar@{.>}[ull]^<<<<<<{} 
&& 1-1+ \ar@{>}[ull]^3 \ar@{.>}[ullll]^<<<<<<<<{} \\
+11- \ar@{>}[ur]^1 \ar@{>}[urrr]_<<<<<<<<3 
&& 11+- \ar@{>}[ul]_<<<<<2 \ar@{>}[urrr]_<<<<<<3  
&& +-11 \ar@{>}[ul]^<<<<<2 \ar@{>}[ur]^1  
&& -+11 \ar@{>}[ur]^<<<<<2 \ar@{>}[ul]^1   
&& 11-+ \ar@{>}[ur]_<<<<<2 \ar@{>}[ulll]_<<<<<<3    
&& -11+ \ar@{>}[ul]^1 \ar@{>}[ulll]^<<<<<<<<3 \\
++--\ar@{>}[u]^2 
 & & +-+- \ar@{>}[ull]_2 \ar@{>}[u]^1\ar@{>}[urr]_<<<<<<3
 && -++- \ar@{>}[ull]^<<<<<<<1\ar@{>}[urr]_<<<<<<3
 && +--+ \ar@{>}[ull]^<<<<<<3\ar@{>}[urr]_<<<<<<1
 && -+-+ \ar@{>}[ull]^<<<<<<3 \ar@{>}[u]^1\ar@{>}[urr]^2
 && --++\ar@{>}[u]^2 
}
}
}
$$
\caption{$\U(2,2)$}
\end{figure}

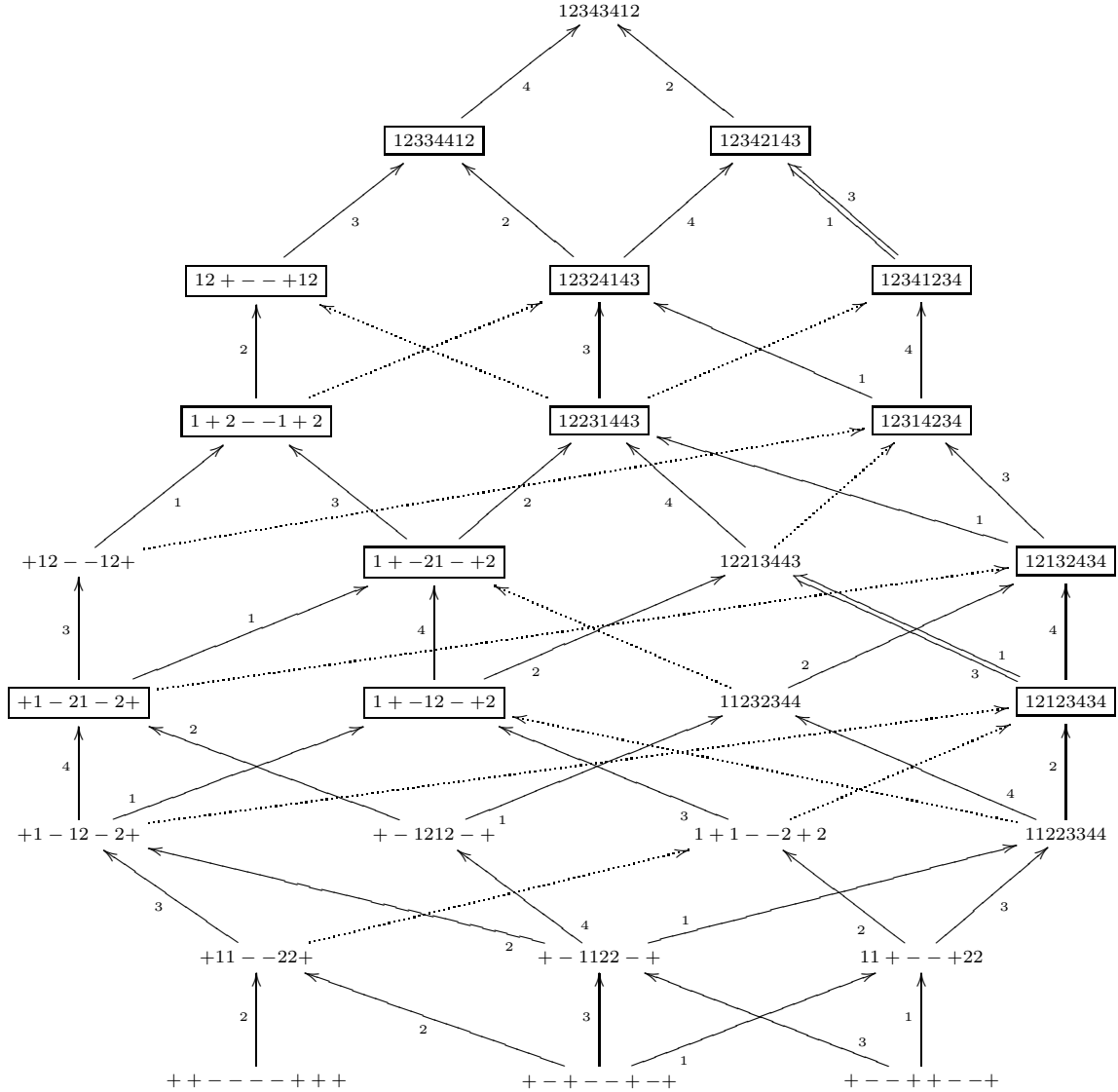
\begin{figure}
\label{f:sp22-example}
$${\tiny
{
\xymatrixcolsep{0pc} 
\xymatrixrowsep{3pc}
\xymatrix
{
&&&12343412\\
&& \boxed{12334412} \ar@{>}_4[ur]&& 
\boxed{12342143}\ar@{>}^2[ul]\\
& \boxed{12+--+12} \ar@{>}_3[ur]
&& \boxed{12324143} \ar@{>}^2[ul] \ar@{>}_4[ur]&& 
\boxed{12341234} \ar@{=>}^1_3[ul]\\
& \boxed{1+2--1+2} \ar@{>}^2[u]
\ar@{.>}_<<<<<<{}[urr]
&& \boxed{12231443} \ar@{>}^3[u] \ar@{.>}_<<<{}[ull]
\ar@{.>}^<<<<{}[urr]
&& 
\boxed{12314234}\ar@{>}^4[u] \ar@{>}_<<<1[ull]\\
+12--12+ \ar@{>}_1[ur]
\ar@{.>}_<<<<<<<{}[urrrrr]
&& \boxed{1+-21-+2} \ar@{>}^3[ul] \ar@{>}_2[ur]
&& 12213443 \ar@{>}^4[ul] \ar@{.>}_<<<<{}[ur]
&& \boxed{12132434} \ar@{>}_3[ul] \ar@{>}_<<<<<1[ulll] \\
\boxed{+1-21-2+} \ar@{>}^3[u] \ar@{>}^1[urr]
\ar@{.>}_<<<<<<<<<<<<<{}[urrrrrr]
&& \boxed{1+-12-+2} \ar@{>}^4[u] \ar@{>}_<<<<<<<2[urr]
&& 11232344 \ar@{>}^<<<<2[urr] \ar@{.>}^<<<<<<{}[ull]
&& \boxed{12123434} \ar@{>}^4[u]
\ar@{=>}_<<<<{1}^<<<<<<{3}[ull] \\
+1-12-2+  \ar@{>}^4[u] \ar@{>}^<<<<1[urr] \ar@{.>}_<<<<<<<<<<{}[urrrrrr]
&& +-1212-+ \ar@{>}_>>>>>>2[ull] \ar@{>}_<<<<<1[urr]
&& 1+1--2+2 \ar@{>}^<<<<<<3[ull] \ar@{.>}_<<<<<{}[urr]&& 
11223344\ar@{>}_<<<<<4[ull] \ar@{>}^2[u] \ar@{.>}^<<<<<<<<<<{}[ullll]\\
&+11--22+  \ar@{>}^3[ul] \ar@{.>}_<<<{}[urrr]
&& +-1122-+ \ar@{>}^<<<<<2[ulll] \ar@{>}_<<4[ul] \ar@{>}^<<<<<<1[urrr]
 && 11+--+22 \ar@{>}^<<<<<<2[ul] \ar@{>}_3[ur]\\
& ++----+++ \ar@{>}^2[u]
&& +-+--+-+ \ar@{>}^2[ull] \ar@{>}^3[u]\ar@{>}_<<<<<<<1[urr]
&& +--++--+ \ar@{>}_<<<<<3[ull] \ar@{>}^{1}[u]
\\
}
}
}
$$
\caption{$\mathrm{PSp}(2,2)$}
\end{figure}

\begin{figure}
\label{f:so*6-example}
$${\tiny
{
\xymatrixcolsep{2pc} 
\xymatrixrowsep{3pc}
\xymatrix@d
{
&&&12+-12 \\
&& 1+21-2  \ar@{>}[ur]^2 && 
1-12+2  \ar@{>}[ul]_3  \\
& +1212-  \ar@{>}[ur]^1 
&& 11-+22  \ar@{>}[ul]_3  \ar@{>}[ur]_2 
&& -1122+   \ar@{>}[ul]_1 \\
+++---  \ar@{>}[ur]^3 
&& 
+--++-  \ar@{>}[ul]_3  \ar@{>}[ur]^1 && 
-+-+-+  \ar@{>}[ul]_1  \ar@{>}[ur]^2  
&& --+-++  \ar@{>}[ul]_2 \\
}
}
}
$$
\caption{$\SO^*(6)$}
\end{figure}

\begin{figure}
\label{f:so*8-example}
{\tiny
{
\xymatrixcolsep{2pc} 
\xymatrixrowsep{3pc}
\xymatrix@d
{
&&&&&&& AABB \\
&&&& \boxed{AA++}  \ar@{>}[urrr]_4
&&& \boxed{ABAB}  \ar@{>}[u]_2&&& \boxed{AA--} \ar@{>}[ulll]_4\\
&& \boxed{A+A+}  \ar@{>}[urr]_2 \ar@{.>}_{}[urrrrr]
&&&& \boxed{ABBA}  \ar@{=>}[ur]^3_1
& & abba  \ar@{>}[ul]_4 \ar@{.>}^<<<<<{}[ullll] \ar@{.>}_<<<<<{}[urr]
&&&& \boxed{A-A-} \ar@{>}[ull]_2
\ar@{.>}_{}[ulllll]
\\
+ AA+ \ar@{>}[urr]^1
\ar@{.>}[urrrrrr]_<<<<<{}
&& A++A \ar@{>}[u]^3 \ar@{.>}[urrrr]_{}
&& \boxed{a+-a} \ar@{>}[ull]^<<<<<4 \ar@{>}[urrrr]^2
&&& \boxed{abab} \ar@{>}[ul]^4 \ar@{=>}[ur]^3_1
&&& \boxed{a-+a} \ar@{>}[ull]^2 \ar@{>}[urr]^4
&& A--A \ar@{>}[u]^3
 \ar@{.>}[ullllll]^<<<<<<{}
&& -AA- \ar@{>}[ull]^1
\ar@{.>}[ullllllll]^<<<<<<<<{}
\\
& +A+A \ar@{>}[ul]^3 \ar@{>}[ur]^<<<<<1
&& +a-a \ar@{>}[ulll]^<<<<<4 \ar@{>}[ur]^<<<<<1
\ar@{.>}[urrrr]_<<<<<{}
&& a+a- \ar@{>}[ulll]^<<<<<4 \ar@{>}[ul]^3
\ar@{.>}[urr]^<<<<<{}
&& aabb \ar@{>}[u]^2 \ar@{.>}^<<<<<{}[ulll]
\ar@{.>}_<<<<<{}[urrr]
&& a-a+ \ar@{>}[ur]^3 \ar@{>}[urrr]_<<<<4
\ar@{.>}[ull]_<<<<<{}
&& -a+a \ar@{>}[ul]_<<<<<1\ar@{>}[urrr]^<<<<<4
\ar@{.>}[ullll]^<<<<<{}
&& -A-A \ar@{>}[ul]^<<<<<1 \ar@{>}[ur]^3 
\\
++AA \ar@{>}[ur]^2
&& +aa- \ar@{>}[ul]^4 \ar@{>}[ur]^3 \ar@{>}[urrr]^<<<<<1
&& +-aa \ar@{>}[ul]^>>>>>2 \ar@{>}[urrr]^<<<<<1
 && aa+- \ar@{>}[ul]^<<<<<2 \ar@{>}[ur]^3
 && aa-+ \ar@{>}[ul]^3 \ar@{>}[ur]^<<<<<2
 && -+aa \ar@{>}[ulll]^<<<<<1 \ar@{>}[ur]_<<<<<2
 && -aa+ \ar@{>}[ulll]^<<<<<1 \ar@{>}[ul]^3 \ar@{>}[ur]^4
 && --AA \ar@{>}[ul]^2 
\\
++++ \ar@{>}[u]^4
& & ++-- 
\ar@{>}[ull]^4 \ar@{>}[u]^2
&& +-+- \ar@{>}[ull]^2 \ar@{>}[u]^3 \ar@{>}[urr]^<<<<<1
&& +--+ \ar@{>}[ull]^<<<<<3 \ar@{>}[urr]_<<<<<1
&& -++- \ar@{>}[ull]^<<<<<1 \ar@{>}[urr]_<<<<<3
&& -+-+ \ar@{>}[ull]^<<<<<1 \ar@{>}[u]^3 \ar@{>}[urr]^2
&& --++ \ar@{>}[u]^2 \ar@{>}[urr]^4
&& ---- \ar@{>}[u]^4\\
}
}
}
\caption{$\SO^*(8)$}
\end{figure}

\medskip

\begin{example}
\label{e:bad}
We conclude with some examples where isogeny considerations necessarily make formulating 
pattern avoidance results more complicated and less uniform.
Let $G = \Sp(2n,\C), \ol G = \PSp(2n,\C) = G/F$, $K = \GL(n,\C)$, and $\ol K = K/F$.   
By considerations similar to those treated in Section
\ref{s:sppq}, orbits of $K$ on $\caB$ are parametrized by $\Sigma^\asym_\pm(2n)$, the union over all $p + q =n$
of antisymmetric elements in $\Spq$.  Orbits of $\ol K$ are parametrized by equivalence classes in $\Sigma_\pm^\asym(2n)$
for the relation generated by $\gamma \sim -\gamma$ and the obvious version of \eqref{e:breaks} holds.  This time,
however, the relationship between the (rational) smoothness of the closure of $\caO_{\bar \gamma}$
and the (rational) smoothness of the closures of its connected components is a little complicated.  For instance,
let $\gamma_2 = 1+-1$, $\gamma_3= +1+-1-$, and $\gamma_4 = 112+-233$.  Then
each $\gamma_i$ contains (in the sense of Definition \ref{d:po}) the pattern $1+-1$.
The $K$ orbits $\caO_{\gamma_i}$ each have closures which are smooth
but not rationally smooth.  Meanwhile the $\ol K$ orbit $\caO_{\ol \gamma_1}$ has closure which
is not rationally smooth, $\caO_{\ol \gamma_2}$ has closure which is rationally smooth (but not
smooth), and $\caO_{\ol \gamma_3}$ has closure which once again is not rationally smooth.
(One may prove the rational smoothness assertions using Theorem \ref{t:springer}; so, in particular,
the criterion of the theorem is sensitive to isogeny.)  Further calculation suggest a relatively simple pattern
avoidance criterion for (rational) smoothness of $K$ orbit closures may exist, but formulating such a result
for $\ol K$ orbit closures is messier.  The situation is similarly complicated for $\SO(p,q)$ (and even more so
if $p+q$ is divisible by four).  This example suggests that it is perhaps reasonable to assume that $G$
is simply connected (and thus $K$ is connected) when formulating pattern avoidance results in general. 
\end{example}

\end{document}